\documentclass[12pt,twoside]{article}
\usepackage[english]{babel}
\usepackage[latin1]{inputenc}
\usepackage{amsmath}
\usepackage{amssymb,amsfonts}
\usepackage{graphicx}
\usepackage{times,amssymb,amscd}
\newtheorem{thm}{Theorem}[section]

\newtheorem{lem}[thm]{Lemma}

\newtheorem{defn}[thm]{Definition}
\newtheorem{rem}[thm]{Remark}
\numberwithin{equation}{section}

\newcommand{\bA}{\mathbf{A}}

\newcommand{\bE}{\mathbf{E}}

\newcommand{\bH}{\mathbf{H}}

\newcommand{\bL}{\mathbf{L}}

\newcommand{\bR}{\mathbf{R}}
\newcommand{\bS}{\mathbf{S}}
\newcommand{\bV}{\mathbf{V}}

\newcommand{\ba}{\mathbf{a}}
\newcommand{\bp}{\mathbf{p}}

\newcommand{\be}{\mathbf{e}}

\newcommand{\bx}{\mathbf{x}}

\newcommand{\bT}{\mathbf{T}}

\newcommand{\bt}{\mathbf{t}}

\newcommand{\BV}{\boldsymbol{V}}

\newcommand{\Be}{\boldsymbol{e}}
\newcommand{\Bu}{\boldsymbol{u}}
\newcommand{\Bv}{\boldsymbol{v}}

\newcommand{\cP}{\mathcal{P}}
\newcommand{\cS}{\mathcal{S}}
\newcommand{\cT}{\mathcal{T}}
\newcommand{\cR}{\mathcal{R}}

\newcommand{\EUC}{\mathbf E^3}
\newcommand{\SPH}{\bS^3}
\newcommand{\HYP}{\bH^3}
\newcommand{\SXR}{\bS^2\!\times\!\bR}
\newcommand{\HXR}{\bH^2\!\times\!\bR}
\newcommand{\SLR}{\widetilde{\bS\bL_2\bR}}
\newcommand{\NIL}{\mathbf{Nil}}
\newcommand{\SOL}{\mathbf{Sol}}

\begin{document}
\pagestyle{myheadings}
\markboth{\centerline{Géza Csima}}
{Isoptic surfaces of segments in $\SXR$ and $\HXR$ geometries}
\title
{Isoptic surfaces of segments in $\SXR$ and $\HXR$ geometries
}

\author{Géza Csima \\
\normalsize Department of Geometry, Institute of Mathematics,\\
\normalsize Budapest University of Technology and Economics, \\
\normalsize M\H{u}egyetem rkp. 3., H-1111 Budapest, Hungary \\
\normalsize csgeza@math.bme.hu
\date{\normalsize{\today}}}

\maketitle
\begin{abstract}

In this work, we examine the isoptic surface of line segments in the $\SXR$ and $\HXR$ geometries, which are from the 8 Thurston geometries. Based on the procedure first described in \cite{CsSz6}, we are able to give the isoptic surface of any segment implicitly. We rely heavily on the calculations published in \cite{Sz202,Sz22}.
As a special case, we examine the Thales sphere in both geometries, which are called {\it Thaloid}. 
In our work we will use the projective model of $\SXR$ and $\HXR$ described by E. Moln\'ar in \cite{M97}.
\end{abstract}

\section{Introduction} \label{section1}

Of the Thurston geometries, those with constant curvature (Euclidean $\EUC$, hyperbolic $\HYP$, spherical $\SPH$)  have been extensively studied, but the other five geometries, $\HXR$, $\SXR$, $\NIL$, $\SLR$, $\SOL$ have been thoroughly studied only from a differential geometry and topological point of view. However, classical concepts can be formulated highlighting the beauty and underlying structure of these geometries, such as: geodesic curves and spheres, equidistant surfaces, translation curves an spheres, lattices, geodesic and translation triangles and their surfaces, their interior angle sum, locus of points from a segment subtends a given angle (isoptic surfaces) and similar statements to those 
known in constant curvature geometries. These question have not been in the focus of attention yet.

In this paper among the 8 Thursten geometries (see \cite{S} and \cite{T}) we are interested in $\SXR$ and $\HXR$ spaces, at the same time. This is primarily due to the fact that a significant symmetry can be observed in both the calculations and the results. That is why in the following sections, after introducing the appropriate notations, we examine the two spaces simultaneously. In Section 2, we will introduce both geometries.

In $\SXR$ and $\HXR$ geometry, we put the question on the agenda of the locus of the points in space from which a given section subtends a given angle. To do this, we need to understand how two points can be connected in these geometries. We basically have three options for this. The Euclidean segment imagined in the model of the corresponding geometries is practical and easy to handle, but has no real geometric  importance in some geometries, in $\SXR$ and $\HXR$ for instance, so we will not deal with it hereafter.

Our second option is the geodesic curve, usually defined as having locally minimal arc length between any two (near enough) points. The equation system of the parametrized geodesic curve $g(x(t),y(t),z(t))$ can be determined by the Levy-Civita theory of Riemann geometry. We can assume, that the starting point of a geodesic curve is the origin because we can transform it into an arbitrary starting point by an appropriate translation. The above procedure gives the geodesic curve as the solution of a second-order differential equation.

The third and last option so far is the translation curve, that in the Thurston spaces, can be introduced in a natural way (see \cite{MoSzi10, Sz12}) by translations mapping each point to any point. Consider a unit vector at the origin. Translations, postulated at the
beginning carry this vector to any point by its tangent mapping. If a curve $t\rightarrow (x(t),y(t),z(t))$ has just the translated
vector as tangent vector in each point, then the  curve is called a {\it translation curve}. This assumption leads to a system of first order 
differential equations, thus translation curves are simpler than geodesics. 

One can ask that is it possible that these curves differ from each other?  The answer is positive generally (excluding some special cases) in $\NIL$, $\SLR$ and $\SOL$ geometries. Moreover, they play an important role and often seem to be more natural in these geometries, than their geodesic lines (see e.g \cite{PSSz1,PSSz2}). In the remaining five Thurston geometries $\EUC,$ $\SPH,$ $\HYP,$ $\SXR$ and $\HXR,$ the translation and geodesic curves coincide with each other. Furthermore, in $\EUC,$ $\SPH$ and $\HYP,$ all three curves are the same and in $\SLR$, translation curves looks like Euclidean straight lines in the model, described in \cite{M97}. From now on, when two points are connected, it is done with the corresponding translation (or geodesic) curve and not with a Euclidean line segment. In Section 2, we recall these curves, determined in \cite{MoSzi10}.

In Section 3, we prepare the matrix of the transformation that pulls an arbitrary point back to the origin of the model and we recall the definition of the isoptic curves and surfaces described in \cite{CsSz6}, for which we use triangles. With this approach, we can avoid problems arising from orientation. Internal angle sum for triangles in $\SXR$ and $\HXR$ had been studied in \cite{Sz202}. We can draw from this study and use the angle calculation method provided there to calculate a single internal angle. This approach seems 
generally effective to determine the implicit equation of the isoptic surface for a segment. We recall the procedure more precisely in Section 3, that have already applied for $\NIL$ in \cite{CsSz6}. We visualize these surfaces in both $\SXR$ and $\HXR$, using their models and Thaloids will be also analyzed as a special case.

\section{On $\SXR$ and $\HXR$ geometries and their translation curves}
In \cite{M06} E.~Moln\'ar has shown that the homogeneous 3-spaces have a unified interpretation in the projective 3-sphere $\mathcal{P} \mathcal{S}^3(\BV^4,\BV_4,\mathbb{R}).$ In this work, we will use this projective model of $\SXR$ and $\HXR.$
We will use the Cartesian homogeneous coordinate simplex $E_0(\be_0)$,$E_1^{\infty}(\be_1)$,$E_2^{\infty}(\be_2)$,
$E_3^{\infty}(\be_3),$ $(\{\be_i\}\subset \bV^4$ with the unit point $E(\be = \be_0 + \be_1 + \be_2 + \be_3 ))$ 
which is distinguished by an origin $E_0$ and by the ideal points of coordinate axes, respectively. 
Moreover, $\mathbf{y}=c\bx$ with $0<c\in \mathbb{R}$ (or $c\in\mathbb{R}\setminus\{0\})$
defines a point $(\bx)=(\mathbf{y})$ of the projective 3-sphere $\mathcal{P} \mathcal{S}^3$ (or that of the projective space $\cP^3$ where opposite rays
$(\bx)$ and $(-\bx)$ are identified). 
The dual system $\{(\Be^i)\}, \ (\{\Be^i\}\subset \BV_4)$, with $\be_i\Be^j=\delta_i^j$ (the Kronecker symbol), describes the simplex planes, especially the plane at infinity 
$(\Be^0)=E_1^{\infty}E_2^{\infty}E_3^{\infty}$, and generally, $\Bv=\Bu\frac{1}{c}$ defines a plane $(\Bu)=(\Bv)$ of $\cP \cS^3$
(or that of $\cP^3$). Thus $0=\bx\Bu=\mathbf{y}\Bv$ defines the incidence of point $(\bx)=(\mathbf{y})$ and plane
$(\Bu)=(\Bv)$, as $(\bx) \text{I} (\Bu)$ also denotes it. Thus $\SXR$ can be visualized in the affine 3-space $\bA^3$
(so in $\bE^3$) and the points of $\HXR$ form an open cone solid, described by the following set:
\begin{equation}
\HXR=\left\{X(\bx=x^i\be_i)\in\mathcal{P}^3:-(x^1)^2+(x^2)^2+(x^3)^2<0<x^0,x^1\right\}
\label{2.1}
\end{equation}

In this context E. Moln\'ar \cite{M97} has derived the well-known infinitesimal arc-length squares invariant under translations at any point  as follows
\begin{equation}
   \begin{gathered}
     \SXR:\ (ds)^2=\dfrac{(dx)^2+(dy)^2+(dz)^2}{x^2+y^2+z^2}\\
     \HXR:\ (ds)^2=\dfrac{1}{(x^2-y^2-z^2)^2}\left\{(x^2+y^2+z^2)(dx)^2+\right.\\
		+2x(dx)(-2y(dy)-2z(dz))+(x^2+y^2-z^2)(dy)^2+\\
		\left.+2yz2(dy)(dz)+(x^2-y^2+z^2)(dz)^2\right\}
       \end{gathered} \label{2.2}
\end{equation}
		
In both geometries, we introduce a new coordinate system in order to write both the arc-length squares and the translation curves more simply. We introduce the polar parametrization of $\SXR$ and the cylindrical parametrization of $\HXR$ in $\BV^4$
\begin{equation}
\begin{gathered}
\SXR:\  x_0=1,\ x_1 = e^t\cos\theta,\ x_2 =e^t \sin\theta\cos\phi,\ x_3 =e^t\sin\theta\sin\phi \\
-\pi<\phi\leq\pi,\ 0\leq\theta\leq\pi,\ t\in\mathbb{R}
\label{2.3}
\end{gathered}
\end{equation}
\begin{equation}
\begin{gathered}
\HXR:\  x_0=1,\ x_1 = e^t\cosh r,\ x_2 =e^t \sinh r\cos\phi,\ x_3 =e^t\sinh r\cos\phi \\
-\pi<\phi\leq\pi,\ 0\leq r,\ t\in\mathbb{R}
\label{2.4}
\end{gathered}
\end{equation}
where $(\theta,\phi)$ and $(r,\phi)$ are the usual polar coordinates of $\bS^2$ and $\bH^2,$ furthermore $t$ is the real component, the so-called fibre coordinate in the direct product $\SXR$ and $\HXR.$ With $x = \frac{x_1}{x_0},$ $y = \frac{x_2}{x_0},$ $z = \frac{x_3}{x_0},$ setting $t$ to be $0$ describes the unit sphere in \ref{2.3}, and the $x>0$ part of the two-sheeted hyperboloid $x^2-y^2-z^2=1$ in \ref{2.4}. This last surface can be called the unit hyperboloid of $\HXR.$ In both geometries $t=\infty$ would be the ideal plane $(\be^0)$ at infinity, $t = -\infty$ would be the origin $(\be_0)$ in limit in $\EUC$ model. Central similarity with factor $e^a$ means the translation by $\mathbb{R}$-component $a,$ commuting with any isometry of $\bS^2$ and $\bH^2.$

Then with the new coordinates, the arc-length square is more simple. 
\begin{equation}
   \begin{gathered}
     \SXR:\ (ds)^2=(dt)^2+(d\theta)^2+\sin^2(\theta)(d\phi)^2\\
     \HXR:\ (ds)^2=(dt)^2+(dr)^2+\sinh^2(r)(d\phi)^2
       \end{gathered} \label{2.5}
\end{equation}

We can assume that the starting point of a translation curve in both geometries is the $(1,1,0,0)$ point, as it is the origin of the coordinate systems, described above. Hereafter, let the functions $S(t)$ and $C(t)$ be $\sin(t)$ and $\cos(t)$ in $\SXR,$ $\sinh(t)$ and $\cosh(t)$ in $\HXR.$ Then the translation curve by \cite{MoSzi10} can be given:
\begin{equation}
\begin{gathered}
x(\tau)=e^{\tau\sin (v)}C(\tau\cos (v)),\\
y(\tau)=e^{\tau\sin (v)}S(\tau\cos (v))\cos(u),\\
z(\tau)=e^{\tau\sin (v)}S(\tau\cos (v))\sin(u),\\
-\pi<u\leq\pi,\ -\frac{\pi}{2}\leq v\leq\frac{\pi}{2}.
\label{2.6}
\end{gathered}
\end{equation}

\begin{figure}[htp]
	\centering
		\includegraphics[width=0.48\textwidth]{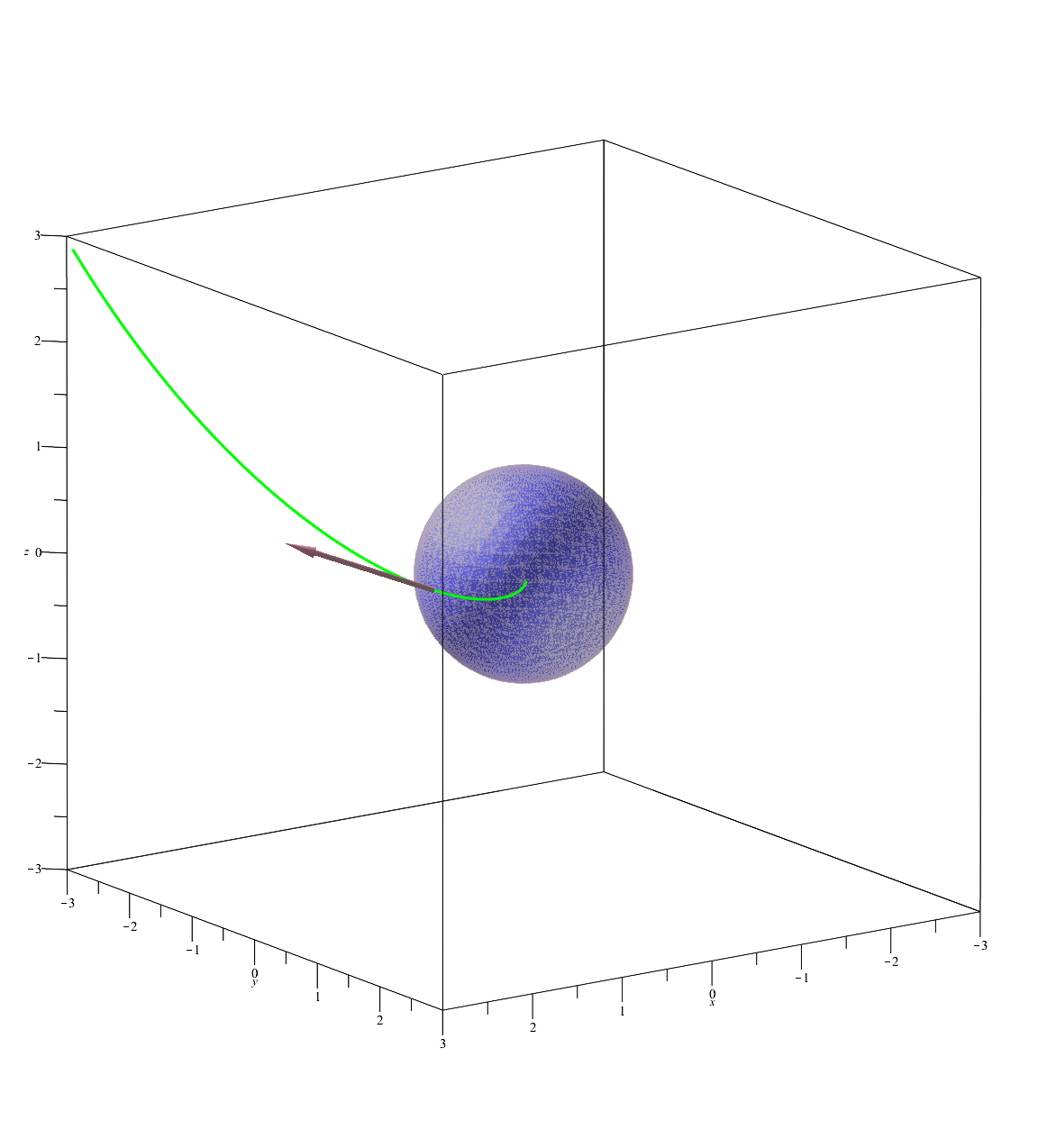} \includegraphics[width=0.48\textwidth]{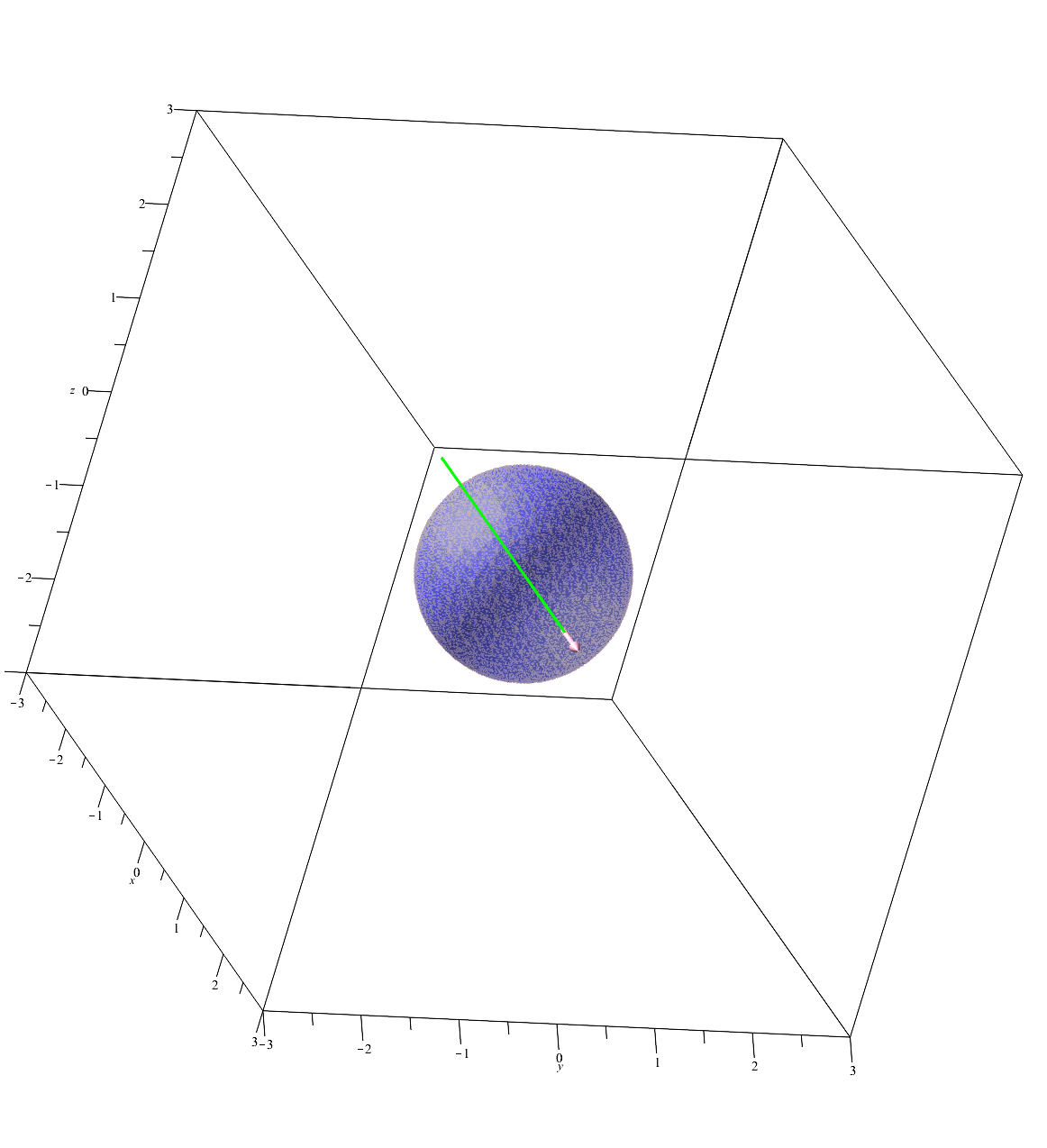}
		\caption{Translation curve in $\SXR$}
	\label{fig:sxr_tr}
\end{figure}

\begin{figure}[htp]
	\centering
		\includegraphics[width=0.48\textwidth]{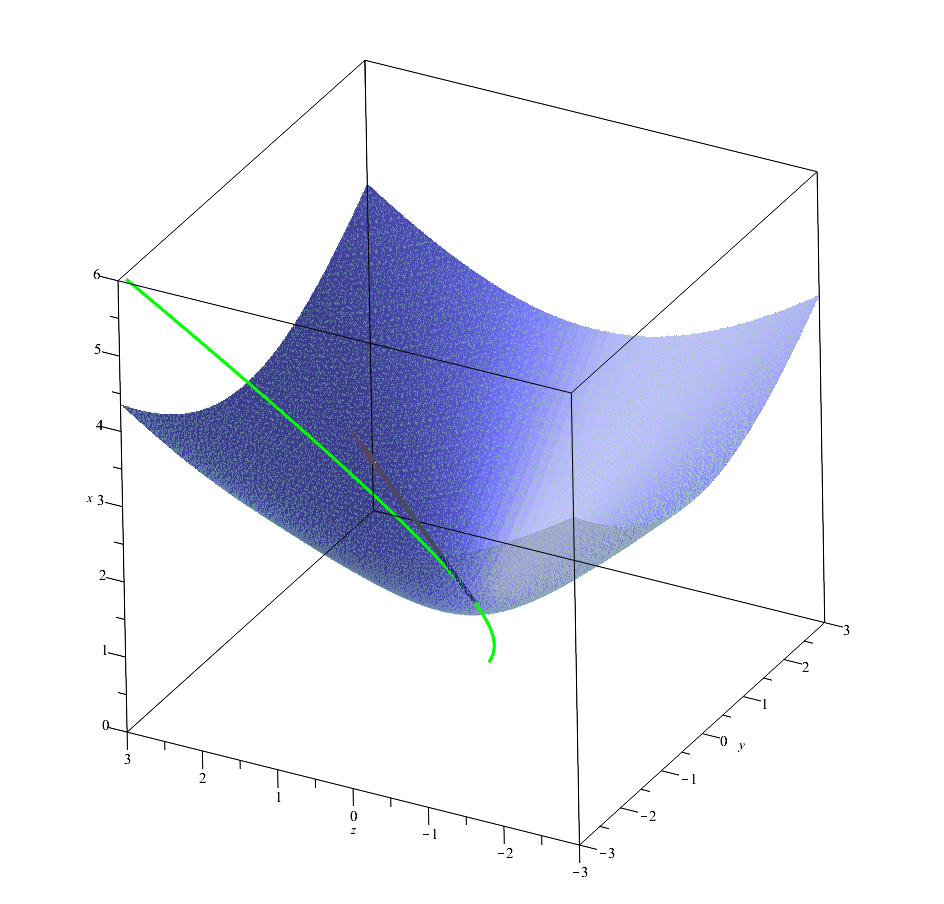} \includegraphics[width=0.48\textwidth]{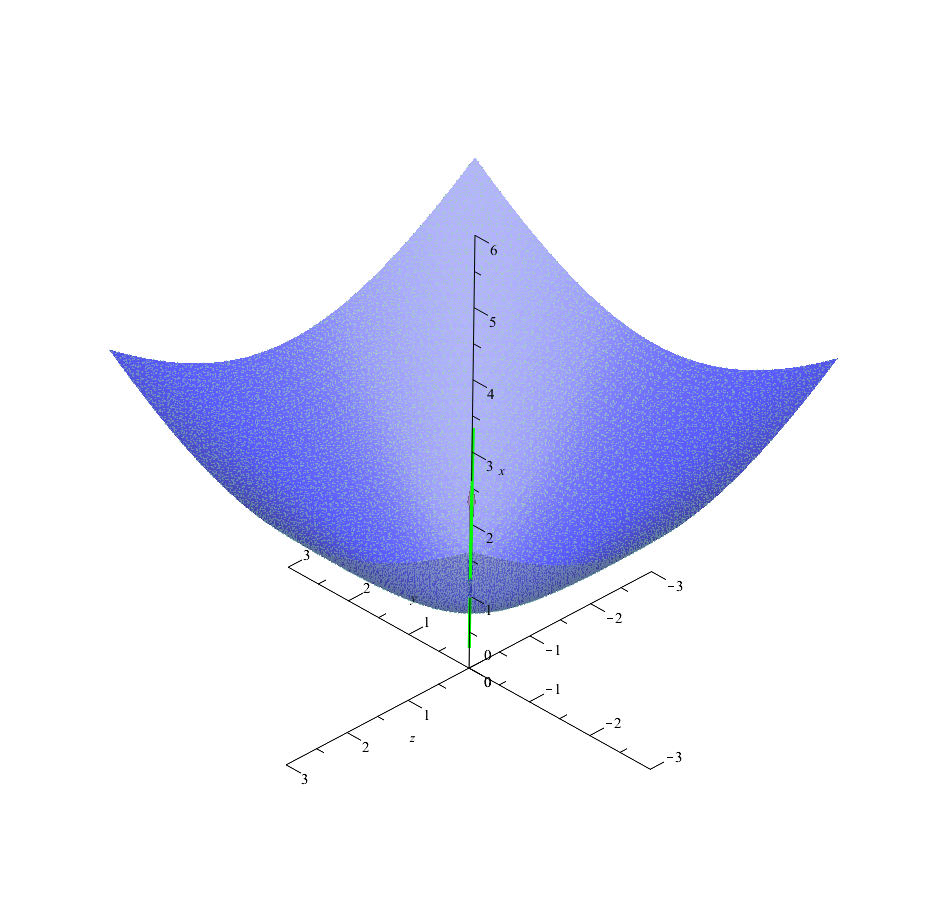}
		\caption{Translation curve in $\HXR$}
	\label{fig:hxr_tr}
\end{figure}

In the parametric equation of the translation curve above $\tau$ denotes the arc-length parameter; $v$ denotes the angle, formed by the starting direction vector of the curve and the tangent plane at the origin $A_1=(1,1,0,0)$ of the unit sphere ($x^2+y^2+z^2=1$) for $\SXR$ and the tangent plane at $A_1$ of unit hyperboloid ($x^2-y^2-z^2=1$) for $\HXR$; and $u$ denotes the angle, formed by the $y$ axis and the projection of this starting direction onto the tangent plane, described above (see the left sides of Figure \ref{fig:sxr_tr} and \ref{fig:hxr_tr}).
\begin{rem}
\begin{enumerate}
	\item It is easy to see, that the translation curve lies in a plane (see the right sides of Figure \ref{fig:sxr_tr} and \ref{fig:hxr_tr}) with equation:
\begin{equation}
\sin(u)y-\cos(u) z=0  
\label{2.7}
\end{equation}
\item The tangent vector of \ref{2.6} at the origin is provided by replacing 0 for $\tau$ in the derivative by $\tau.$
\begin{equation}
\bt=(\sin(v),\cos(v)\cos(u),\cos(v)\sin(u)) 
\label{2.8}
\end{equation}
\end{enumerate}
\end{rem}


\begin{defn}
The distance $d(P_1,P_2)$ between the points $P_1$ and $P_2$ is defined by the arc length of the above (see \ref{2.6}) translation curve from $P_1$ to $P_2$.
\label{2.9}
\end{defn}

\begin{defn} The sphere of radius $R >0$ with centre at the origin, (denoted by $S_O(R)$), with the usual longitude and altitude parameters 
$u$ and $v$, respectively by \ref{2.6}, is specified by the following equations:
\begin{equation}
\begin{gathered}
        S_O(R): \left\{ \begin{array}{ll} 
       x(\tau)=e^{R\sin v}C(R\cos v),\\
			 y(\tau)=e^{R\sin v}S(R\cos v)\cos(u),\\
			 z(\tau)=e^{R\sin v}S(R\cos v)\sin(u).
        \end{array} \right.
        \label{2.10}
\end{gathered}
\end{equation}
\end{defn}

\begin{defn}
 The body of the sphere of centre $O$ and of radius $R$ in the $\SXR$ and $\HXR$ spaces is called ball, denoted by $B_{O}(R)$,
 i.e. $Q \in B_{O}(R)$ iff $0 \leq d(O,Q) \leq R$.
\end{defn}

The parametrization in \ref{2.10} allows us, to create the implicit equation of $S_{O}(R)$:
 \begin{equation}
\begin{gathered}
        \dfrac{\log^2(x^2\pm(y^2+z^2))}{4}+\textrm{arcC}^2\left(\dfrac{x}{x^2\pm(y^2+z^2)}\right)=R^2,
						\label{2.11}
\end{gathered}
\end{equation}
where $\pm$ and $\textrm{arcC}(x)$ is $+$ and $\arccos(x)$ for $\SXR,$ $-$ and $\mathrm{arccosh}(x)$ for $\HXR.$

\section{Isoptic surface}

\subsection{Transformation matrix}

In the following, we describe a series of transformations from the isometry group of $\SXR$ and $\HXR$ such that their composition translate an arbitrary point $P(1,x,y,z)$ back to the origin of the model, i.e. $E_0(1,1,0,0)$ in both geometries. This is important for us because, as can be seen from the arc-length square formula (see \ref{2.5}), at this point the angles appear in their real size, and we also chose this point as the starting point of each translation curve.

First, we perform a so called fibre translation $\cT,$ that shifts the point onto the unit sphere in $\SXR$ and onto the unit hyperboloid in $\HXR.$ We remind the dear reader that $\pm$ is $+$ for $\SXR$ and $-$ for $\HXR.$

\begin{equation}
\cT=diag\left\{1,\dfrac{1}{\sqrt{x^2\pm(y^2+z^2)}},\dfrac{1}{\sqrt{x^2\pm(y^2+z^2)}},\dfrac{1}{\sqrt{x^2\pm(y^2+z^2)}}\right\}
\label{3.1}
\end{equation}

\begin{equation}
P^\cT=\left(1,\dfrac{x}{\sqrt{x^2\pm(y^2+z^2)}},\dfrac{y}{\sqrt{x^2\pm(y^2+z^2)}},\dfrac{z}{\sqrt{x^2\pm(y^2+z^2)}}\right)
\label{3.2}
\end{equation}

Then we rotate around the $x$ axis so that the last coordinate of the image is zero:

\begin{equation}
\cR_x=
\begin{pmatrix}
1 & 0 & 0 & 0 \\
0 & 1 & 0 & 0 \\
0 & 0 & \frac{y}{\sqrt{y^2+z^2}} & -\frac{z}{\sqrt{y^2+z^2}} \\
0 & 0 &  \frac{z}{\sqrt{y^2+z^2}} & \frac{y}{\sqrt{y^2+z^2}} \\
\end{pmatrix}
\label{3.3}
\end{equation}

\begin{equation}
P^{\cT\cR_x}=\left(1,\dfrac{x}{\sqrt{x^2\pm(y^2+z^2)}},\dfrac{\sqrt{y^2+z^2}}{\sqrt{x^2\pm(y^2+z^2)}},0\right)
\label{3.4}
\end{equation}

It can be verified that the plane of the translation curve drawn from $E_0$ to $P$ changes under $\cR_x$ and transforms into the plane $z=0$ (see \ref{2.7}). Then we rotate it around the $z$ axis so that the second to last coordinate of the image is also zero. During this transformation, the plane of the translation curve does not change.
\begin{equation}
\cR_z=
\begin{pmatrix}
1 & 0 & 0 & 0 \\
0 & \frac{x}{\sqrt{x^2\pm(y^2+z^2)}} & -\frac{\sqrt{y^2+z^2}}{\sqrt{x^2\pm(y^2+z^2)}} & 0 \\
0 & \pm\frac{\sqrt{y^2+z^2}}{\sqrt{x^2\pm(y^2+z^2)}} & \frac{x}{\sqrt{x^2\pm(y^2+z^2)}} & 0 \\
0 & 0 & 0 & 1
\end{pmatrix}
\label{3.5}
\end{equation}

\begin{equation}
P^{\cT\cR_x\cR_z}=\left(1,1,0,0\right)
\label{3.6}
\end{equation}

Finally, to make the plane of the transformed and the original translation curve coincide, we apply the inverse of the $\cR_x$ transformation.
\small
\begin{equation}
\begin{gathered}
\bT=\cT\cR_x\cR_z\cR_x^{-1}=\\
\begin{pmatrix}
1 & 0 & 0 & 0 \\
0 & \dfrac{x}{x^2\pm(y^2+z^2)} & -\dfrac{y}{x^2\pm(y^2+z^2)} & -\dfrac{z}{x^2\pm(y^2+z^2)} \\
0 & \pm\dfrac{y}{x^2\pm(y^2+z^2)} & \dfrac{xy^2+z^2\sqrt{x^2\pm(y^2+z^2)}}{(x^2\pm(y^2+z^2))(y^2+z^2)} & \dfrac{xyz-yz^2\sqrt{x^2\pm(y^2+z^2)}}{(x^2\pm(y^2+z^2))(y^2+z^2)} \\
0 & \pm\dfrac{z}{x^2\pm(y^2+z^2)} &  \dfrac{xyz-yz^2\sqrt{x^2\pm(y^2+z^2)}}{(x^2\pm(y^2+z^2))(y^2+z^2)} & \dfrac{xz^2+y^2\sqrt{x^2\pm(y^2+z^2)}}{(x^2\pm(y^2+z^2))(y^2+z^2)} 
\end{pmatrix}
\label{3.7}
\end{gathered}
\end{equation}
\normalsize

\subsection{Isoptic surfaces}

It is well known that in the Euclidean plane the locus of points from a segment subtends a given angle $\alpha$ $(0<\alpha<\pi)$ is the union of two arcs except for the endpoints with the segment as common chord. If this $\alpha$ is equal to $\frac{\pi}{2}$ then we get the Thales circle. Replacing the segment to another general curve, we obtain the Euclidean definition of isoptic curve:
\begin{defn}[\cite{yates}]The locus of the intersection of tangents to a curve meeting at a constant angle $\alpha$ $(0<\alpha <\pi)$ is the $\alpha$ -- isoptic of the given curve. The isoptic curve with right angle called \textit{orthoptic curve}.
\label{defiso}
\end{defn}
\begin{rem}Sometimes we consider the $\alpha$ -- and $\pi-\alpha$ -- isoptics together. Thus, in the case of the section, we get two circles with the segment as a common chord (endpoints of the segment are excluded). Hereafter, we call them $\alpha$ -- isoptic circles.
\end{rem}

Although the name "isoptic curve" was suggested by Taylor  in 1884 (\cite{T}), reference to former results can be found in \cite{yates}. In the obscure history of isoptic curves, we can find the names of la Hire (cycloids 1704) and Chasles (conics and epitrochoids 1837) among the contributors of the subject. A very interesting table of isoptic and orthoptic curves is introduced in \cite{yates}, unfortunately without any exact reference of its source. However, recent works are available on the topic, which shows its timeliness. In \cite{harom} and \cite{tizenketto}, the Euclidean isoptic curves of closed strictly convex curves are studied using their support function.
Papers \cite{Kurusa,Wu71-1,Wu71-2} deal with Euclidean curves having a
circle or an ellipse for an isoptic curve. Further curves appearing as isoptic curves are well studied in Euclidean plane geometry $\mathbf{E}^2$, see e.g. \cite{Loria,Wi}.
Isoptic curves of conic sections have been studied in \cite{H} and \cite{Si}. There are results for Bezier curves by Kunkli et al. as well, see \cite{Kunkli}. Many papers focus on the properties of isoptics, e.g. \cite{nyolc,MM3,het}, and the references  therein. There are some generalizations of the isoptics as well \textit{e.g.} equioptic curves in \cite{Odehnal} by Odehnal or secantopics in \cite{tizenegy, tiz} by Skrzypiec.

We can extend the very first question to the space: "What is the locus of points where a given segment subtends a given angle?" 
Or a question equivalent to the former: "For the given spatial points $A$ and $B$, what is the locus of the points $P$ for which 
the internal angle at $P$ of the triangle $ABP\triangle$ is a given angle?" We use this to define the $\alpha$ -- isoptic surface of a Euclidean spatial segment.

\begin{defn} The $\alpha$ -- isoptic surface of a Euclidean spatial segment $\overline{A_1A_2}$ is the locus of points $P$ for which the internal angle 
at $P$ in the triangle, formed by $A_1,$ $A_2$ and $P$ is $\alpha.$ If $\alpha$ is the right angle, then it is called the Thaloid of $\overline{A_1A_2}.$
\label{defisoszakaszE}
\end{defn}

It is easy to see in the Euclidean space that:
\begin{thm}
The locus of points in the Euclidean space from where a given segment subtends a given angle $\alpha$ $(0<\alpha <\pi)$ or $\pi-\alpha$ 
is a self-intersecting torus obtained by rotating the $\alpha$ -- isoptic circles drawn in any plane containing the section around the line of the section. $\square$
\end{thm}
\begin{rem} 
\begin{enumerate}
\item The torus in the above theorem contains both the isoptic surface for the given angle and the supplementary angle. In this case, we can easily separate the $\alpha$ -- and $\pi-\alpha$ -- isoptic surfaces along the self-intersection.
Specifically, the orthoptic surface is a sphere whose diameter is the section. We can call this the Thaloid of the segment.
\item There is no point in examining the isoptic surface defined in the above way for other spatial curves, because if the curve is not of constant 
0 curvature, then there is an external point from which the curve and the point cannot be fitted into a plane. In this case, the above definition needs to be generalized.
\end{enumerate}
\end{rem}
For further isoptic surfaces in Euclidean geometry, see \cite{CsSz4,CsSz5}, where we extend the definition of isoptic surfaces to other spatial objects.
The notion of isoptic curve can be extended to the other planes of constant curvature (hyperbolic plane $\mathbf{H}^2$ and spherical plane $\mathbf{S}^2$). We studied these questions in \cite{Cs-Sz14} and \cite{Cs-Sz20}. 

\begin{defn} The $\SXR$ or $\HXR$ $\alpha$ -- isoptic surface of a segment $\overline{A_1A_2}$ is the locus of points $P$ 
for which the internal angle at $P$ in the triangle, formed by $A_1,$ $A_2$ and $P$ is $\alpha.$ If $\alpha$ is the right angle, 
then it is called the Thaloid of $\overline{A_1A_2}.$
\label{defisoszakasz}
\end{defn}

We emphasize here that the section itself does not appear in our calculations, we only deal with the endpoints. We can assume by the homogeneity of the geometries that one of its endpoints coincide with the origin $A_1=E_0=(1,0,0,0)$ and the other is $A_2=(1,a,b,c)$. Considering a point $P=(1,x,y,z)$, we can determine the angle $A_1PA_2\angle$ along the procedure described below.

First, we apply $\bT_{P}$ to all three points (see \ref{3.7}). This transformation preserves the angle $A_1PA_2\angle$ and pulls back $P$ to the origin, hence the angle in question seems in real size. 

\begin{equation}
\begin{gathered}
\bT_{P}(P)=(1,1,0,0); \\
\bT_{P}(A_1)=\left(1,\dfrac{x}{x^2\pm(y^2+z^2)},-\dfrac{y}{x^2\pm(y^2+z^2)},-\dfrac{z}{x^2\pm(y^2+z^2)}\right);\\
\bT_{P}(A_2)=\begin{pmatrix} 1\\
\dfrac{ax\pm(by+cz)}{x^2\pm(y^2+z^2)}\\
\dfrac{zy(cx-az)-y^2(ay-bx)+z(bz-cy)\sqrt{x^2\pm(y^2+z^2)}}{(x^2\pm(y^2+z^2))(y^2+z^2)} \\ 
\dfrac{z^2(cx-az)-yz(ay-bx)-y(bz-cy)\sqrt{x^2\pm(y^2+z^2)}}{(x^2\pm(y^2+z^2))(y^2+z^2)} 
\end{pmatrix}^T
\label{eltoltak}
\end{gathered}
\end{equation}

To determine the internal angle in $A_1A_2P_\triangle$ at $P$, we need the tangent vectors of the translation curves running into the points $\bT_{P}(A_1)$ and $\bT_{P}(A_2)$ from $\bT_{P}(P)=E_0.$ It is not necessary to determine the exact value of the parameters $u,$ $v,$ $\tau$, it is enough to evaluate the vector $\bt$ (see \ref{2.8}). 

\begin{lem} 
Let $(1,x,y,z)$ $(x,y,z\in\mathbb{R}$ and $x^2+y^2+z^2>0$ in $\SXR$; $x^2-y^2-z^2>0,\ x>0$ in $\HXR)$ be the homogeneous coordinates of a point $P\in\SXR$ or $P\in\HXR.$ Then the translation curve, drawn to $P$ from $E_0=(1,1,0,0)$ has the following tangent:
\footnotesize
\begin{equation}
\tau\cdot\bt_P=\left(\dfrac{1}{2}\ln(x^2\pm(y^2+z^2)),
\dfrac{y\,\mathrm{arcC}\left(\dfrac{x}{\sqrt{x^2\pm(y^2+z^2)}}\right)}{\sqrt{y^2+z^2}},
\dfrac{z\,\mathrm{arcC}\left(\dfrac{x}{\sqrt{x^2\pm(y^2+z^2)}}\right)}{\sqrt{y^2+z^2}}\right)
\end{equation}
\normalsize
where $\pm$ and $\mathrm{arcC}(x)$ is $+$ and $\arccos(x)$ for $\SXR,$ $-$ and $\mathrm{arccosh}(x)$ for $\HXR,$ and $\tau$ is the distance of $P$ and $E_0.$ $\square$ 
\label{tan_lem}
\end{lem}

Applying the above tangent formula to the points $\bT_{P}(A_1)$ and $\bT_{P}(A_2)$, we obtain the vectors, forming the interior angle at $P.$ 
%
%
Finally, we are able to determine the angle, using the usual angle formula, derived from the dot product of two vectors. Due to the length of the formula, the result is only presented in its seriously simplified form (using appropriate dot and cross products and the spherical/hyperbolic law of cosine for sides), introducing some notations in the following theorem. 
\begin{figure}[htp]
	\centering
		\includegraphics[width=0.58\textwidth]{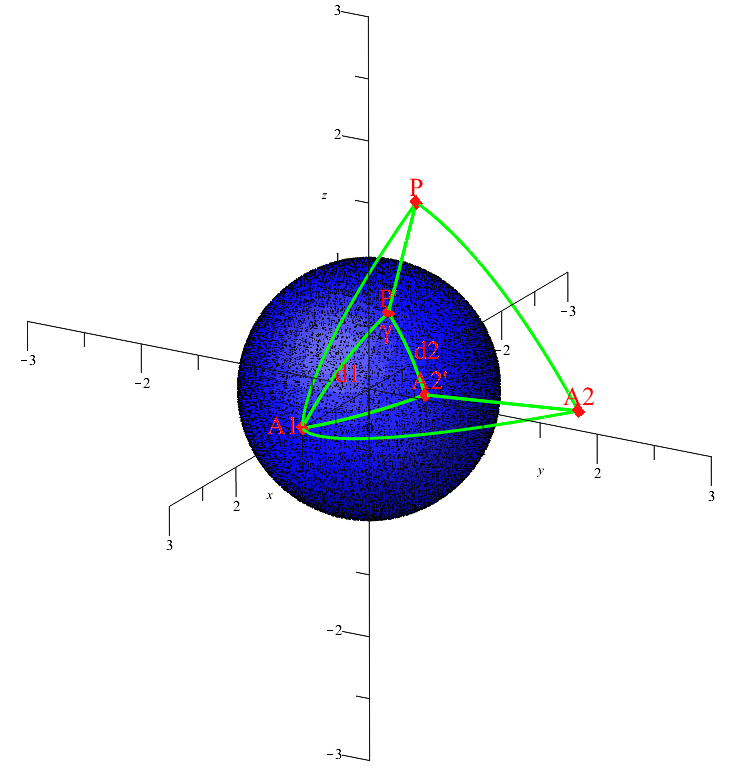} 
		\caption{Explanatory figure for Theorem \ref{isothm} in $\SXR$ with $A_1(1,1,0,0),$ $A_2(2,3,1)$ and $P(1,1,2)$}
	\label{fig:theorem}
\end{figure}

\begin{thm}
Let $A_1=(1,1,0,0)$ and $A_2=(1,a,b,c)$ be given points in $\SXR$ or $\HXR,$ and $\alpha$ be a given angle. If $\ba_1=(1,0,0);$ $\ba_2=(a,b,c);$ $\bp=(x,y,z);$ the projected image of $A_1,$ $A_2$ and $P$ onto the unit surface (sphere in $\SXR$ and hyperboloid in $\HXR$) of the geometry along the corresponding fibre lines are $A_1'=A_1$, $A_2'$  
and $P';$ 
$d_1$ and $d_2$ are the distances of $P'$ to $A_1'$ and to $A_2'$ and $\gamma$ is the internal angles at $P'$ in $A_1'A_2'P'_\triangle$ (see \ref{fig:theorem} in $\SXR$); then the $\alpha$ -- isoptic surface of the $\overline{A_1A_2}$ segment has the equation
\begin{equation}
\cos(\alpha)=\dfrac{d_1d_2\cos(\gamma)+\ln\frac{|\ba_1|}{|\bp|}\cdot\ln\frac{|\ba_2|}{|\bp|}}
{\sqrt{\left(\ln^2\frac{|\ba_1|}{|\bp|}+d_1^2\right)\left(\ln^2\frac{|\ba_2|}{|\bp|}+d_2^2\right)}}.\square
\label{eq:isopt}
\end{equation} 
\label{isothm}
\end{thm}

\begin{figure}[htp]
	\centering
		\includegraphics[width=0.48\textwidth]{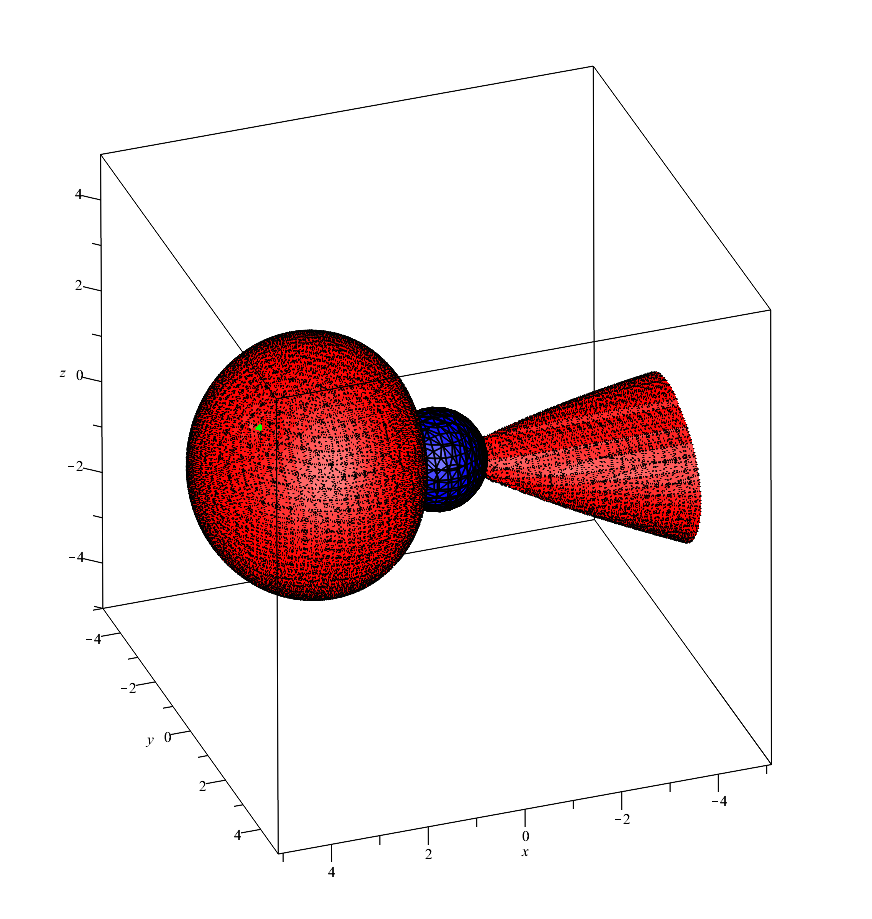} \includegraphics[width=0.48\textwidth]{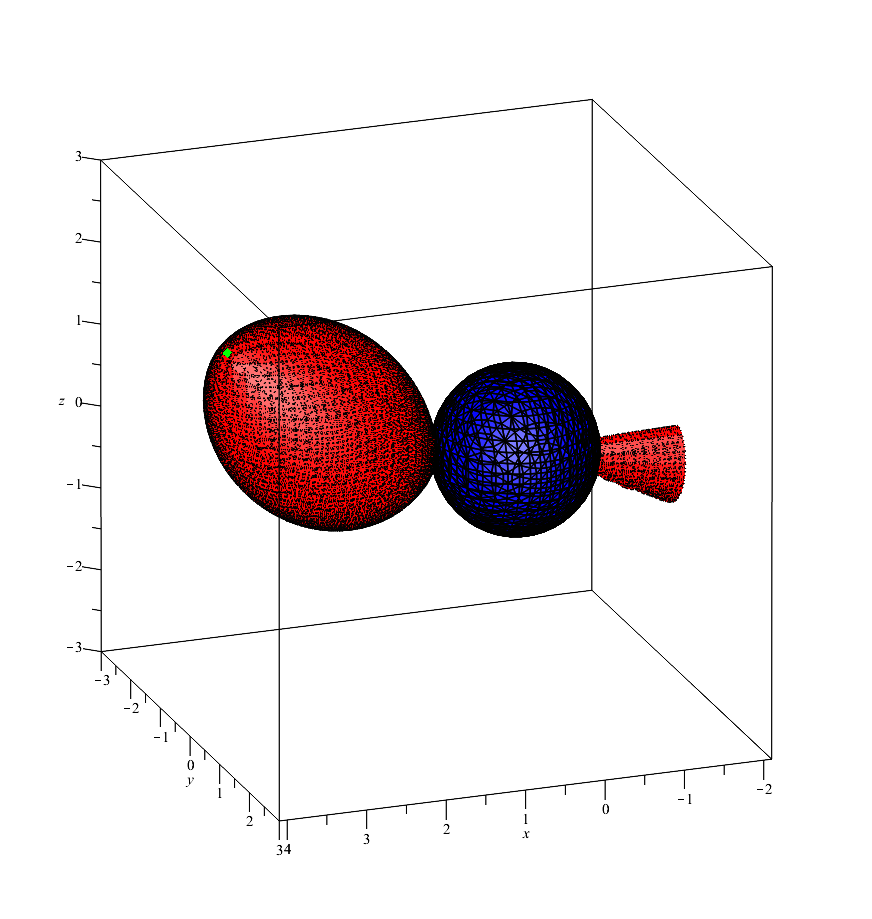}
		\caption{Isoptic surfaces of the $A_1(1,1,0,0)$ and $A_2(1,4,1,2)$ segment in $\SXR$ geometry with $\alpha=80^\circ$ (left) and $\alpha=120^\circ$ (right)}
	\label{fig:sxr}
\end{figure}
\begin{figure}[htp]
	\centering
		\includegraphics[width=0.48\textwidth]{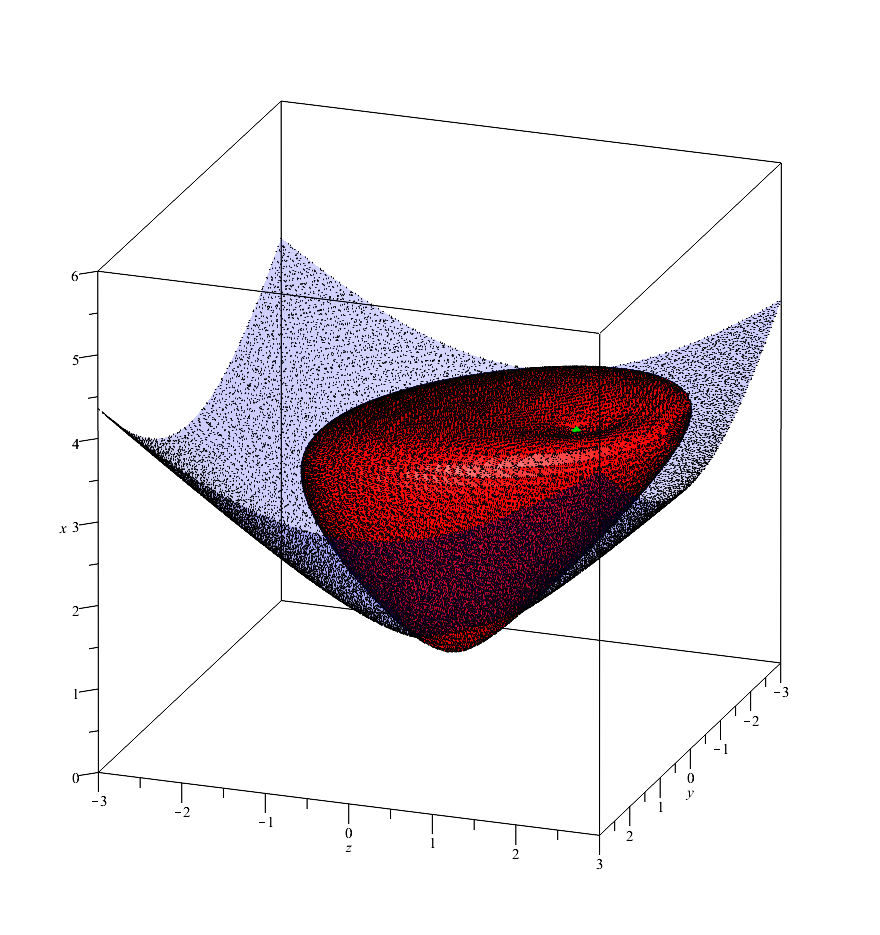} \includegraphics[width=0.48\textwidth]{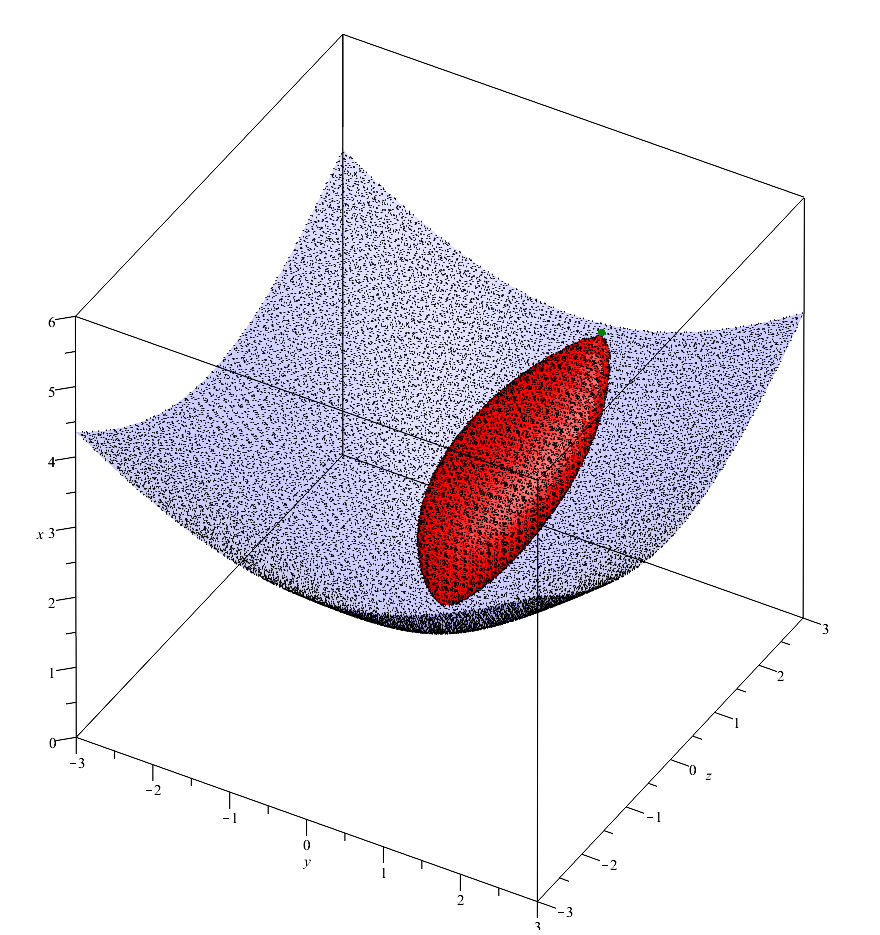}
		\caption{Isoptic surfaces of the $A_1(1,1,0,0)$ and $A_2(1,4,1,2)$ segment in $\HXR$ geometry with $\alpha=75^\circ$ (left) and $\alpha=120^\circ$ (right)}
	\label{fig:hxr}
\end{figure}

\begin{rem} If also $A_2$ lies on the unit surface of the geometry, i.e. $|a_2|=1,$ then restricting $P$ to the unit surface $(|\bp|=1)$ will result in $\bS^2$ or $\bH^2$ geometry and the $\cos(\alpha)=\cos(\gamma)$ equation. As it has mentioned before, $\cos(\gamma)$ was the result of a cosine theorem so that it can be computed just from the sine and cosine (hyperbolic) functions of distances $d_1,$ $d_2$ and $d_0=\mathrm{dist}(A_1;A_2).$ 
\end{rem}

Let us examine the special case when the endpoints of the segment are situated on the $x$ axis, i.e. $A_1=(1,0,0,0)$ and $A_2=(1,a,0,0).$ 
In this case, the segment is along a fibre line and it looks like a Euclidean segment in the model. 
Applying $\bT$ to $A_1$ and $A_2,$ we get that:
\begin{equation}
\begin{gathered}
\bT_P(A_1)=\left(1,\dfrac{x}{x^2+y^2+z^2},-\dfrac{y}{x^2+y^2+z^2},-\dfrac{z}{x^2+y^2+z^2}\right),\\ 
\bT_P(A_2)=\left(1,\dfrac{ax}{x^2+y^2+z^2},-\dfrac{ay}{x^2+y^2+z^2},-\dfrac{az}{x^2+y^2+z^2}\right) 
\label{3.11}
\end{gathered}
\end{equation}

According to Lemma \ref{tan_lem}, we can determine the $\bt_1$ and $\bt_2$ tangents of the translation curves, drawn to $\bT_P(A_1)$ and $\bT_P(A_2)$

\footnotesize
\begin{equation}
\begin{gathered}
\bt_1=\left(\dfrac{1}{2}\ln\left(\dfrac{1}{x^2\pm(y^2+z^2)}\right),
-\dfrac{y\,\mathrm{arcC}\left(\dfrac{x}{\sqrt{x^2\pm(y^2+z^2)}}\right)}{\sqrt{y^2+z^2}},
-\dfrac{z\,\mathrm{arcC}\left(\dfrac{x}{\sqrt{x^2\pm(y^2+z^2)}}\right)}{\sqrt{y^2+z^2}}\right),\\ 
\bt_2=\left(\dfrac{1}{2}\ln\left(\dfrac{a^2}{x^2\pm(y^2+z^2)}\right),
-\dfrac{y\,\mathrm{arcC}\left(\dfrac{x}{\sqrt{x^2\pm(y^2+z^2)}}\right)}{\sqrt{y^2+z^2}},
-\dfrac{z\,\mathrm{arcC}\left(\dfrac{x}{\sqrt{x^2\pm(y^2+z^2)}}\right)}{\sqrt{y^2+z^2}}\right)
\label{3.12}
\end{gathered}
\end{equation}
\normalsize

Since we are interesting in the Thaloid, where the angle of $\bt_1$ and $\bt_2$ is $\frac{\pi}{2},$ we consider their dot product to be zero:
\begin{equation}
\left\langle\bt_1,\bt_2 \right\rangle=\dfrac{1}{4}\ln\left(\dfrac{1}{x^2\pm(y^2+z^2)}\right)\ln\left(\dfrac{a^2}{x^2\pm(y^2+z^2)}\right)+
\mathrm{arcC}^2\left(\dfrac{x}{\sqrt{x^2\pm(y^2+z^2)}}\right)=0
\label{3.13}
\end{equation}

To better understand the nature of the above implicit surface, let us apply a translation which pulls back the midpoint of $\overline{A_1A_2}$ to $A_1.$ The translation curve to $A_2$ has a very simple parametrization in this case: $(e^{t},0,0),$ where $t\in(0,\ln(a)).$ Then the coordinates of the midpoint is $F=(1,\sqrt{a},0,0).$ The appropriate fibre translation,that maps $A_1$ to $F$ is 
$\cT=\mathrm{diag}\left\{1,\sqrt{a},\sqrt{a},\sqrt{a}\right\},$ so that $x=x'\sqrt{a},$ $y=y'\sqrt{a}$ and $z=z'\sqrt{a}.$ Then \ref{3.13} has a different form:
\begin{equation}
\begin{gathered}
\dfrac{1}{4}\ln\left(a((x')^2\pm((y')^2+(z')^2))\right)\ln\left(\dfrac{(x')^2\pm((y')^2+(z')^2)}{a}\right)+\\
+\mathrm{arcC}^2\left(\dfrac{x'}{\sqrt{(x')^2\pm((y')^2+(z')^2)}}\right)=0 \ \ \Longleftrightarrow \\
\dfrac{1}{4}\ln^2\left((x')^2\pm((y')^2+(z')^2)\right)+\mathrm{arcC}^2\left(\dfrac{x'}{\sqrt{(x')^2\pm((y')^2+(z')^2)}}\right)=\dfrac{1}{4}\ln^2(a)
\label{3.14}
\end{gathered}
\end{equation}

Summarizing the results above, we get that:

\begin{lem}
Let $A_1=(1,1,0,0)$ and $A_2=(1,a,0,0)$ $(a\in\mathbb{R}^+)$ be given points in $\SXR$ or $\HXR.$ Then the Thaloid of the $\overline{A_1A_2}$ segment is a sphere with centre $C=(1,\sqrt{a},0,0)$ and radius $r=\ln(\sqrt{a}).$
\end{lem}

\begin{figure}[htp]
	\centering
		\includegraphics[width=0.48\textwidth]{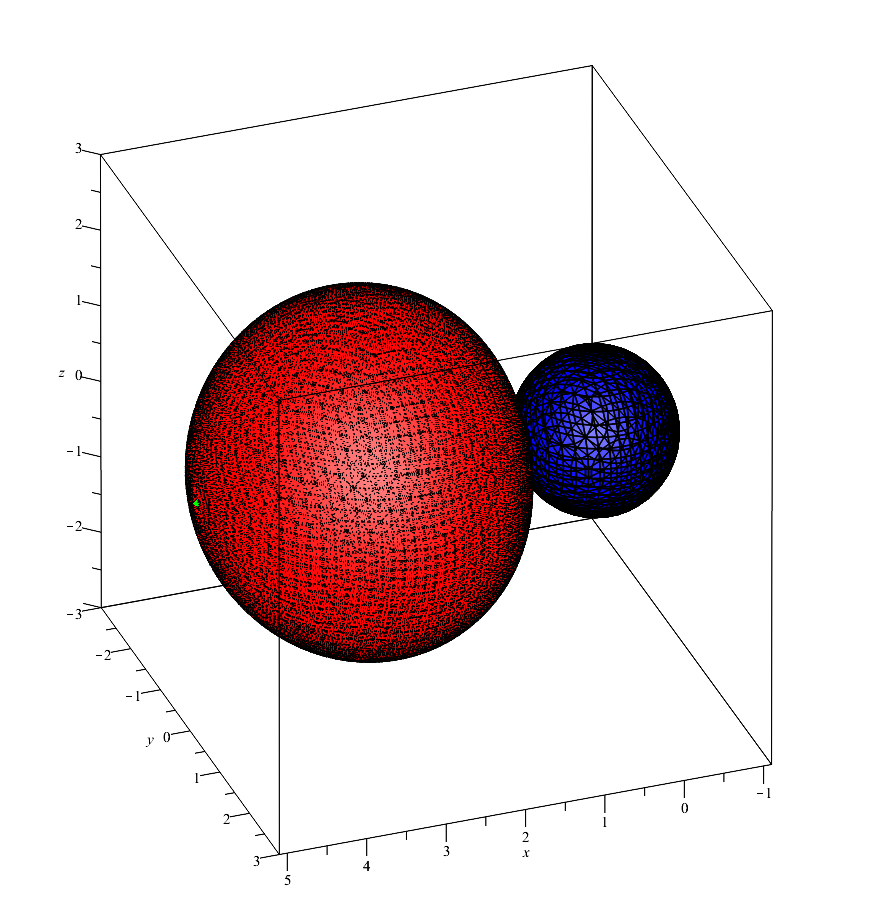} \includegraphics[width=0.48\textwidth]{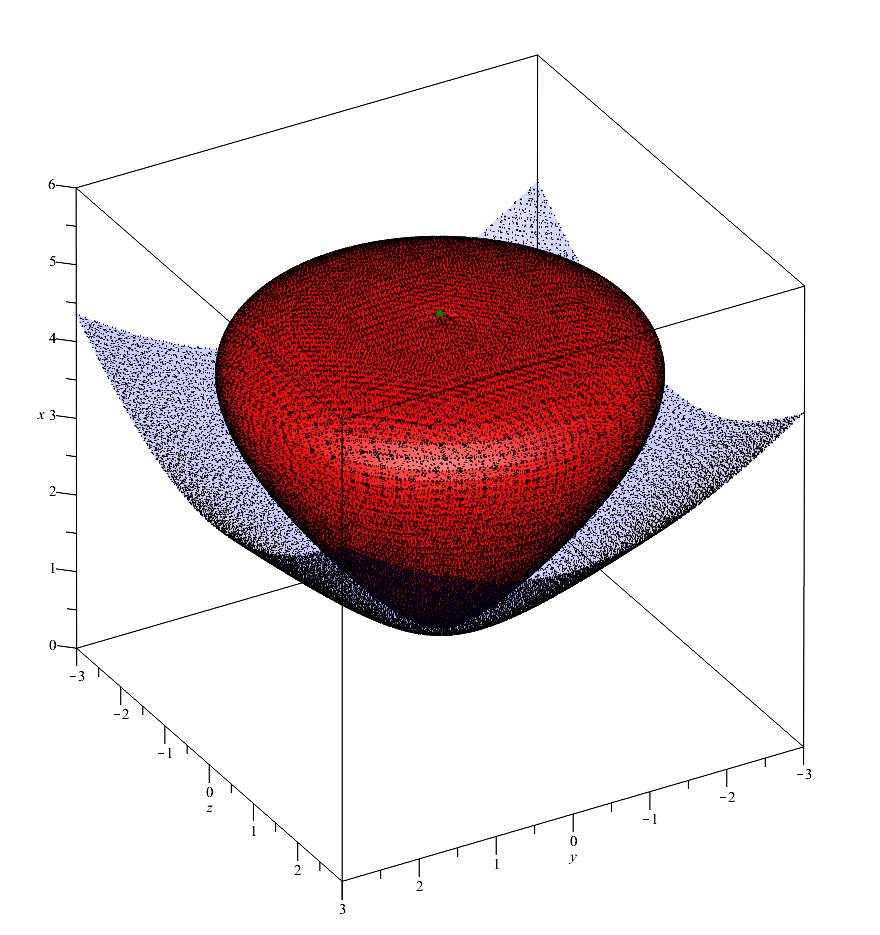}
		\caption{Thaloid of the $A_1(1,1,0,0)$ and $A_2(1,5,0,0)$ segment in $\SXR$ (left) and $\HXR$ (right) geometries}
	\label{fig:Thaloid}
\end{figure}

\section*{Data Availability Statement}
Data sharing not applicable to this article as no databases were generated or analyzed during the current study.

\section*{Conflict Of Interest Statement}
The authors have no conflict of interest to declare that are relevant to the content of this study.

\end{document}